%
%
\input amstex
\documentstyle{gen-p}
\NoBlackBoxes

\issueinfo {00}
{}
  {}
 {XXXX}

\topmatter
\title A sequence of connections and a characterization of
K\"ahler manifolds \endtitle
\author Mikhail Shubin\endauthor
\leftheadtext{MIKHAIL SHUBIN}%
\rightheadtext{Connections and K\"ahler manifolds}

\address Department of Mathematics, Northeastern University, Boston,
MA 02115 \endaddress


\email shubin\@neu.edu\endemail

\thanks The author was partially supported 
by NSF grant DMS-9706038. \endthanks


\dedicatory Dedicated to Mel Rothenberg 
on the occasion of his 65th birthday.\enddedicatory


\abstract
We study a sequence of connections which 
is associated with a Riemannian metric 
and an almost symplectic structure on a manifold $M$ 
according to a construction in [5].  We prove that 
if this sequence is trivial (i.e. constant)  or 
2-periodic, then $M$ has a canonical K\"ahler structure.
\endabstract

\endtopmatter

\document

\def\sqr#1#2{{\vcenter{\vbox{\hrule height.#2pt
    \hbox{\vrule width.#2pt height#1pt \kern#1pt
    \vrule width.#2pt}
    \hrule height.#2pt}}}}

\def\Ga{\Gamma}

\def\om{\omega}

\def\pa{\partial}

\font\bbb=msbm10

\def\C{\hbox{{\bbb C}}}

\def\ms{\medskip}

\head Introduction \endhead

In this paper we consider a sequence of connections which is
associated with two structures on the same manifold $M$: 
a Riemannian metric $g$ and an almost symplectic structure $\om$. This sequence was 
introduced in [5]. We discuss a special situation when this sequence is in fact
trivial, i.e. contains only one connection. We prove that in this case
$M$ is a K\"ahler manifold and the K\"ahler structure is obtained from
the structures $(g,\om)$ by a simple procedure which includes a Gromov
construction of an almost complex structure $J$ which 
is compatible with $\om$. 
In our case $J$ proves to be integrable and together with $\om$ (which
proves to be closed) they generate a K\"ahler structure.

\ms
We also note that if the sequence of connections has period 2, then
it is in fact trivial, so we again have a K\"ahler structure. 

\ms
In  Sections 1 and 2 we recall the construction of the sequence 
of connections from [5] and the Gromov construction of an almost complex
structure associated with a pair $(g,\om)$.  The characterization of the 
K\"ahler manifolds is given in Sect.3.

\ms
It is well known that the existence of a K\"ahler structure imposes
many topological restrictions on $M$. In particular, 
there are symplectic manifolds $(M,\om)$ which do not admit
any  K\"ahler structure, even with a different $\om$ 
(see e.g.  [2], p.147, for a Thurston example, and also 
the book [9] for many other examples). 

\ms
It would be interesting to 
provide topological obstructions to the existence of a pair of forms
$(g,\om)$ such that the corresponding sequence of connections is
periodic with a specific period. The simplest new situations arise when
this sequence is trivial but the initial connection is not symmetric 
or when the sequence (starting with the Levi-Civita connection of $g$)
has period 2 not from the very beginning but from some further term.

\ms
All periodic situations give rise to new classes of manifolds which are
generalizations of K\"ahler manifolds.  It is quite probable that 
they are of geometric interest.

\head 1. Sequence of connections associated with two forms.
\endhead

\medskip
Let $M$ be a $C^\infty$ manifold, $\dim M=n$,  
$g$ is a pseudo-riemannian structure on $M$, i.e.
a non-degenerate symmetric form on each tangent space $T_xM$, $x\in M$, which 
smoothly depends on $x$. 
In local coordinates $(x^1,\dots, x^n)$ such a form is given by 
a non-degenerate symmetric tensor $g_{ij}=g(\pa_i,\pa_j)$, where $\pa_j=\pa/\pa x_j$ 
is a basic
vector field in the domain $U$ of the local coordinates. Here $g_{ij}=g_{ij}(x)$
is a $C^\infty$ function on $U$. 

\medskip
Let $\nabla$ be a linear connection i.e. a connection in the tangent bundle
$TM$. It is given locally by the Christoffel symbols $\Ga^k_{ij}$ defined 
from the relation $\nabla_{\pa i}\pa_j=\Ga^k_{ij}\pa_k$, where we use 
the standard summation convention. We will identify the connection $\nabla$ 
with the set of its Christoffel symbols, and will talk about a connection $\Ga$
(instead of $\nabla$) which does not lead to a confusion.

\medskip
The torsion tensor of the connection $\nabla$ is given by 
$T^k_{ij}=\Ga^k_{ij}-\Ga^k_{ji}$.

\medskip
The connection $\nabla$ preserves $g$ if 
$$
Xg(Y,Z)=g(\nabla_X Y,Z)+g(Y,\nabla_X Z), \tag 1.1
$$
where $X,Y,Z$ are arbitrary vector fields on $M$.

\medskip
It is well known due to T.~Levi-Civita, H.~Weyl and E.~Cartan (see [3]) 
that for any  tensor $T=(T^k_{ij})$, which is 
antisymmetric with respect to the subscripts $i,j$, there exists a unique connection
$\Ga=(\Ga^k_{ij})$ with the torsion tensor $T$ (i.e. $\Ga^k_{ij}-\Ga^k_{ji}=T^k_{ij}$ ) 
such that $\Ga$ preserves $g$. If $T=0$, this connection is called 
the {\it Levi-Civita connection}. It is the {\it canonical} connection associated
with $g$. 

\ms
A skew-symmetric (or rather almost symplectic) analog of this fact was 
established in [5]. Namely, let $\om$ be an almost symplectic structure on $M$
i.e. a  non-degenerate 2-form on $M$ or skew-symmetric bilinear form on each
tangent space $T_xM$ which is smooth in $x$. In local coordinates such a form is 
determined by a non-degenerate skew-symmetric tensor $\om_{ij}=\om(\pa_i,\pa_j)$.
A connection $\nabla$ (or $\Ga$) on $TM$ preserves 
$\om$ if an identity similar to (1.1)
holds, namely:
$$
X\om(Y,Z)=\om(\nabla_X Y,Z)+\om(Y,\nabla_X Z), \tag 1.2
$$
for any vector fields $X,Y,Z$ on $M$.

\ms
Now for any connection $\Gamma=(\Gamma_{ij}^k)$ on $M$ denote by $\Pi=(\Pi_{ij}^k)$ its
symmetric part ($\Pi=\Pi(\Gamma)$) i.e.
$$\Pi_{ij}^k={1\over 2}(\Gamma_{ij}^k+\Gamma_{ji}^k).\tag 1.3$$
\medskip
\proclaim{Theorem 1.1} {\rm [5]} The map 
$\Gamma\mapsto\Pi(\Gamma)$ gives a bijective
affine correspondence between the set of all connections 
preserving $\omega$ and 
the set of all symmetric connections. The inverse map $\Pi\mapsto\Gamma$ is
given by
$$\Gamma_{kij}={1\over 2}(\pa_k\omega_{ij}-\pa_i\omega_{jk}-\pa_j\omega_{ki})+
(\Pi_{kij}+\Pi_{jik}-\Pi_{ijk}).\tag 1.4$$    
In case when $\omega$ is closed (hence symplectic) this formula can be
rewritten as
$$\Gamma_{kij}=\pa_k\omega_{ij}+(\Pi_{kij}+\Pi_{jik}-\Pi_{ijk}). 
\tag 1.5$$
\endproclaim

\remark{Remark 1.2} It is a well known fact that a symmetric connection preserving 
$\omega$ exists if and only
if $\omega$ is closed (see e.g. [7], [8], [10], [5]). 
\endremark

\medskip
Now we are ready to recall the construction of the sequence of
connections from [5].

Suppose we are given non-degenerate symmetric and antisymmetric
2-forms $g$ and $\omega$ at the same time. Our first (or rather $0$th) connection 
$\nabla_0$ (or $\Gamma_0$) will be just the 
Levi-Civita connection of $g$. It is symmetric 
(or its torsion is $T_0=0$), so we will also
put $\Pi_0=\Ga_0$. The next connection $\nabla_1$ (or $\Gamma_1$)
will be the unique connection preserving $\om$ and having 
the symmetric part $\Pi(\Ga_1)=\Pi_0$.  Denote its torsion $T_1$.
The next connection will be the unique connection $\nabla_2$
(or $\Gamma_2$) which preserves $g$ and has the torsion $T_1$.
Denote by $\Pi_1$ its symmetric part. Then we can take
the connection $\nabla_3$ (or $\Gamma_3$) which preserves $\om$
and has the symmetric part $\Pi_1$. In this way we obtain a sequence of connections
$$
\Ga_0,\Ga_1,\Ga_2,\dots, \tag 1.6
$$
which is uniquely defined by the following properties:

\ms
(i) The connection $\Ga_{2k}$ preserves $g$ for any $k=0,1,\dots$.

(ii) The connection $\Ga_{2k+1}$ preserves $\om$ for any $k=0,1,\dots$.

(iii) The connections $\Ga_{2k}$ and $\Ga_{2k+1}$ have the same
symmetric part for any $k=0,1,\dots$.

(iv) The connections $\Ga_{2k+1}$ and $\Ga_{2k+2}$ have the same 
torsion tensor for any $k=0,1,\dots$.

(v) The connection $\Ga_0$ is symmetric.

\ms
Instead of requiring (v) 
we could also start with an arbitrary torsion-type tensor $T_0$
and take $\Ga_0$ preserving $g$ and having the torsion tensor $T_0$.
This gives us a sequence which also depends on the choice of $T_0$.

\ms
Note that for K\"ahler manifolds, starting with $T_0=0$ we obtain the constant sequence
i.e. we get $\Ga_0=\Ga_1=\Ga_2=\dots$.  Our main goal here will be to
prove the inverse statement under the condition that 
the  metric $g$ is in fact positive (i.e. Riemannian).

\head 2. Gromov's almost complex structure construction
\endhead

\ms
Let $(M,\om)$ be an almost symplectic manifold. Let us recall 

\definition{Definition 2.1} An almost complex structure $J$ in $TM$ is called
{\it compatible with the almost symplectic structure} $\om$, if 

(i) $\om(X,JX)>0$ for any tangent vector $X\not= 0$.

(ii) $\om(JX,Y)+\om(X,JY)=0$ for any tangent vectors $X,Y\in T_xM$, $x\in M$,

\noindent
or, equivalently: 
$$
g(X,Y)=\om(X,JY) \leqno (2.1)
$$
 is an hermitian metric i.e. a 
Riemannian  metric such that $J$ is an isometry in this metric.
\enddefinition

Note also that a condition equivalent to (ii) is
$$
\om(JX,JY)=\om(X,Y).
$$

\proclaim
{Proposition 2.1} {\rm (M.~Gromov [4])}  On any almost symplectic manifold 
$(M,\om)$  there exists  an almost complex structure $J$ which is compatible with $\om$.
\endproclaim

\demo{Proof}
In fact M.~Gromov gave a construction of  such an almost complex structure 
$J$ in $TM$. This construction depends 
on a choice of an arbitrary  start-up Riemanian metric $g_0$. After this choice
the construction of $J$ becomes canonical. So in fact we need a linear algebra
construction in each fiber of $TM$ which is canonical in $\om$, $g_0$
(in particular smoothly depending on $\om$ and $g_0$). Let us describe this
construction following [1] (section 6.1.1).

Denote $V=T_xM$ where $x\in M$ is fixed, so $V$ is a symplectic vector space.
Define a linear endomorphism $A:V\to V$ by
$$
g_0(AX,Y)=\om(X,Y).
$$
Then $A$ is antisymmetric with respect to $g_0$. Therefore $-A^2$ is symmetric and 
$$
g_0(-A^2X,X)=-\om(AX,X)=\om(X,AX)=g_0(AX,AX)>0,
$$
provided $X\in V\setminus 0$. Therefore $-A^2$ is a symmetric positive definite
operator. Define $B=\sqrt{-A^2}$ to be the symmetric positive-definite square-root of
$-A^2$. This can be done canonically e.g. by using the Cauchy integral and choosing 
a branch of the square root in $\C\setminus (-\infty,0]$ with $\sqrt{1}=1$.
Clearly $B$ is invertible and commutes with $A$. Now take $J=B^{-1}A$. Then 
$J^2=-\hbox{Id}_V$. It is also clear that $J$ is antisymmetric 
with respect to $g_0$. Furthermore,
$$
\om(X,JX)=g_0(AX,B^{-1}AX)>0,\quad X\in V\setminus 0,
$$
because $B^{-1}$ is also symmetric and positive definite. 
We also have
$$
\om(JX,Y)=g_0(AJX,Y)=g_0(AB^{-1}AX,Y)
=-g_0(B^{-1}AX,AY)
$$
$$
=-g_0(AX,B^{-1}AY)=-\om(X,JY),
$$
which proves that $J$ is compatible with $\om$.

It follows also that the formula (2.1) defines an hermitian 
metric $g$ on $M$.
\quad\qed
\enddemo

\head 3. Main results
\endhead

Let us assume that we have a triple $(M,g_0,\om)$, where $M$ is a 
 $C^\infty$-manifold, $g_0$ is a Riemannian metric on $M$, 
$\om$ is an almost-symplectic form on $M$. Let $\nabla$ 
(or $\Ga$) denote
the Levi-Civita connection associated with $g_0$. 

\proclaim{Theorem 3.1} Let us assume that the sequence 
of connections on $M$ associated with $M,g_0,\om$ as described in
Sect.1, is trivial, i.e. $\nabla$ preserves $\om$ 
(as well as $g_0$). Then there exists an integrable complex structure $J$ on
$M$ (i.e. on $TM$), such that if $g$ is defined by (2.1), then
$g$ is an hermitian metric on $M$ with 
the same Levi-Civita connection $\nabla$, and the triple 
$(g,\om,J)$ defines a K\"ahler structure on $M$.
\endproclaim

\demo{Proof} Note first that if $\nabla\om=0$ for a symmetric connection 
$\nabla$, then $d\om=0$ (see e.g. [7], [8], [10] and Remark 1.4 in [5]), 
so the form $\om$ is closed, hence it is a symplectic form.

Now take the almost complex structure $J$ and 
the Riemannian metric $g$ which are constructed by the data $(M,\om,g_0)$ as in the proof
of Proposition 2.1. By definition the connection $\nabla$ preserves both $\om$ and $g_0$
(i.e. $\nabla\om=\nabla g_0=0$).
We claim that the same is true for $J$ and $g$.

Indeed, $J$ and $g$
are constructed canonically from $\om$ and $g_0$. 
But the $\nabla$-parallel 
transport (which preserves both $\om$ and $g_0$ ) preserves 
all the relations which are described in the proof of
Proposition 2.1. Hence it preserves $J$ and $g$. In
particular, $\nabla$ is the Levi-Civita connection for $g$
(as well as for $g_0$).

Now we need to establish that $J$ is integrable 
(i.e. $J$ is a complex structure on $M$) and 
$(g,\om,J)$ define  a K\"ahler structure on $M$. This follows from
the arguments given in [6], vol.2, Ch.IX, Sect.4 or [2], p.148.
\quad\qed
\enddemo

\proclaim{Theorem 3.2}
Assume that the sequence of connections associated with
$(g_0,\om)$ has period 2, i.e. $\Ga_2=\Ga_0$. Then in fact
we also have $\Ga_0=\Ga_1=\Ga_2=\dots$, i.e. 
the sequence of connections is trivial and $M$ has a K\"ahler
structure which is canonically defined by $(g_0,\om)$.
\endproclaim

\demo{Proof} The relation $\Ga_2=\Ga_0$ implies in particular
that $\Ga_1$ is symmetric, because $\Ga_1$ and $\Ga_2$ share
the same torsion tensor. On the other hand $\Ga_0$ and $\Ga_1$
share the same symmetric part, so it follows that $\Ga_0=\Ga_1$
because both $\Ga_0$ and $\Ga_1$ are symmetric. 
\quad\qed
\enddemo 

\bigskip
{\it Acknowledgements.} I am grateful to I.~Bernstein, 
S.~Gelfand and B.~Kostant for useful remarks 
and discussions.

\refstyle{A}

\enddocument